%% file: Rama_bis.tex
\input Definitions.tex
\centerline{\bf GENERALIZED PELL-FERMAT EQUATIONS AND PASCAL TRIANGLE}
\medskip
\centerline{Daniel GANDOLFO, Michel ROULEUX}
\smallskip
\centerline{Aix-Marseille Univ, Universit\'e de Toulon, CNRS, CPT, Marseille, France}
\smallskip
\centerline{dgandolfo@wanadoo.fr, rouleux@univ-tln.fr}
\medskip
\noindent{\it Abstract}: Using Pascal triangle, we give a simple generalization to a well-known problem of S.Rama- nujan. 
Thus we are interested in computing the median of
some integer valued distributions, arising naturally when extending partial sums of the arithmetic progression (triangular numbers) to
tetrahedral numbers and beyond. We show this reduces to 
equations of Pell-Fermat type of higher order, which admit very few integer solutions, but 
for which, following S.Ramanujan's original idea, we can always find integer sequences of best approximation, in the Diophantine sense.
In absence of a general theory on Pell-Fermat equation of higher order, our procedure relies much on formal Calculus with Mathematica.
\medskip 
\noindent {\bf 1/ Introduction}
\smallskip
\noindent {\it a) Pascal triangle and progressions}
\smallskip
Pascal triangle identifies with an infinite lower triangular matrix with elements $c_{nk}={n\choose k}$,
$0\leq k\leq n$, the binomial coefficients. It contains many remarkable sequences~: the sum of binomial coefficients in $n$:th row equals $2^n$, 
the sum of diagonal elements are the terms of Fibonacci sequence, and the columns $C_k$ generalize arithmetic sequences.
Namely, the first column $C_0$ is the constant sequence equal to 1~; the second one $C_1$, which gives the partial sums of $C_0$
consists in the arithmetic sequence with ratio 1 and first term 1~; the third one $C_2$ gives the partial sums of $C_1$ by using
the relation ${n\choose 1}+{n\choose2}={n+1\choose2}$~; the fourth one $C_3$ gives the partial sums of $C_2$ by
${n\choose 2}+{n\choose3}={n+1\choose3}$, etc\dots All sequences $C_k$ have a polynomial growth ${\cal O}(n^{k})$.
Because of their simple geometric interpretation,
numbers in $C_2$ are known as the {\it triangular numbers}, these in $C_3$ the {\it tetraedric numbers}, and those occuring in
higher order columns $C_k$ the $k$-{\it simplicial numbers}. 

Pascal triangle also admits asymmetric generalizations, which amount to replace 1 by arbitrary numbers on the second diagonal.
Asymmetric Pascal triangle are built up the same way as standard Pascal triangle, and also extend to complex values of $n$ 
(and complex values of $k$ as well by using Euler $\Gamma$-function).
The simplest one consists in replacing 1 by 
$a\in{\bf C}$, so that the binomial coefficient $c_{nk}$ reads
$$c_{nk}^a={n\choose k}_a={n\choose k}+a{n-1\choose k-1}$$
In particular in the second column $C_1^a$ we recognize the  arithmetic sequence of ratio 1 and first term $1+a$. 
The term $a{n-1\choose k-1}$ accounts for a ``lower order term'' as $n$ becomes large.
\medskip
\noindent {\it b) The ``Houses of Ramanujan''}
\smallskip
This problem, leading to a Pell-Fermat equation was reportledly discovered in quaint circumstances by S.Ramanujan [Ran] and 
solved (exactly) by the method of continued fractions. The history of this Diophantine equation goes to Archimedes, then proceeds to
Bhaskara, Brahmagupta, Fermat, Wallis, Euler and Lagrange who almost brought the theory to its definite form.  

The problem of the ``Houses of Ramanujan'' deals with column $C_2$ of triangular numbers. We say that $m\in{\bf N}^*=\{1,2,3,\cdots\}$
is a {\it House of Ramanujan of order} 2 (or simply a House of Ramanujan) iff the sum ${m\choose2}$ of the 
$m-1$ first integers ``to the left of the House of Ramanujan'' 
is equal to the sum of the next $n-m$ integers ``to the right of the House of Ramanujan'', for some $n> m$, with the 
convention that $n=1$ when $m=1$, that is
$${n+1\choose2}=2{m\choose2}+m\leqno(C2)$$
In the language of Statistics, $m$ is the median of the cumulants $\Sum_{i\leq j}i$, $j\leq n$.
Geometrically this means that for some values $m< n$,
the number of integer points in the triangle $x_1+x_2\leq n$, $x_1,x_2\geq0$ is twice the
number of integer points in the triangle $x_1+x_2\leq m$, $x_1,x_2\geq0$, plus $m$.

We can consider the first quartile as well, defined as the integer $m$ such that
$${n+1\choose2}=4{m\choose2}+m\leqno(Q1)$$
or the third quartile, defined as the integer $m$ such that
$${n+1\choose2}={4\over3}{m\choose2}+m\leqno(Q3)$$
The point of course is to find integer solutions $(n,m)$. The quadratic case $C_2$ plays a peculiar role, and all solutions of (C2)
are given by a sequence in ${\bf Z}[\sqrt2]$ with constant coefficients. 
\smallskip
\noindent {\it c) Main results}
\smallskip
The situation becomes drastically different when trying to find the House of Ramanujan of order $k\geq3$,
i.e. to decompose simplicial numbers of higher order $k\geq3$.  
For $k=3$ this amounts to solve
$${n+2\choose3}=2{m+1\choose3}+{m+1\choose2}\leqno(C3)$$
If this equation has any solution, then the number of integer points in the tetrahedron (with positive coordinates)
$x_1+x_2+x_3\leq n$, 
would be twice the
number of integer points in the smaller tetrahedron $x_1+x_2+x_3\leq m$, plus the triangular number ${m+1\choose2}$. 
Our guess is that there are no solution to (C3), (C4) and (C5). We found only one solution to (C6), corresponding to $(m,n)=(10,11)$. 

The corresponding equations, after some affine change of coordinates, take the form 
$P_k(x,y)=0$, where $P_k$ is polynomial  with integer coefficients (see formulae
(A.2) below), 
i.e. $P_k(x,y)=x^k-Dy^k+\cdots$ (or permuting $x$ with $y$) where the dots mean 
a polynomial of degree $\leq k-2$ without mixed terms, and containing only monomials with the same parity as $k$.

But it turns out that we can always compute asymptotic 
solutions in the Diophantine sense, generalizing Ramanujan's approach. 

Next we address the problem of finding (exact) integer solutions, or pairs $(m,n)$ with half-integer $m$.
We also look at the first and third quartile of such statistical distributions.

In case $C_2$ the method relies on known results for Pell-Fermat equation see [Wo] and references therein. While $u^2-Dv^2=1$, $D$ not a perfect 
square,  has 
always integer solutions given by a (unique) sequence in ${\bf Z}[\sqrt D]$, this holds for $u^2-Dv^2=-1$ only for some 
values of $D$, using continued fractions. Once we have solved $u^2-Dv^2=\pm1$, we can consider the more general Eq. $x^2-Dy^2=c$.
Whenever we guess at a particular solution, we can build a sequence of solutions (called a fundamental sequence),
using the sequences for $u^2-Dv^2=\pm1$. So there exists so many fundamental sequences as ``fundamental particular solutions'' we can 
find for $x^2-Dy^2=c$. But such a family of solutions does not necessarily generate all solutions.

Applying this method for $C_2$ we find all solutions of (C2) (which was of course known before), and one fundamental sequence for (Q3).
Other (related) solutions could be also found by considering the action of some discrete groups acting on the 
hyperbola of equation $x^2-Dy^2=c$, as is the case in the problem of  
quasi-isoceles triangles with a square angle. But we could not find new solutions this way.

Since $P_4$ is a polynomial in $x^2$ and $y^2$, we can reduce $C_4$ to a Pell-Fermat Eq., then we are left to 
find the solutions which
are perfect squares (quadrature). Within the range of values we have 
considered, i.e. $n\leq10^8$, we have shown this way there are no integer solutions 
(or solutions with half-integer $m$), when using the fundamental sequences given by Mathematica. 

We call equation $P_k(x,y)=0$ a {\it generalized Pell-Fermat equation}. Actually nothing seems to be known about these equations,
and since elementary technics, as the reduction of $P_k(x,y)=0$ 
modulo prime numbers $p$ have not led us to any significant issue (see however Sect.3)
we tried to resort to formal calculus on Mathematica. 

So the case $C_3$ and $C_6$ rely instead on the resolution of a 3:rd degree polynomial, which makes use of real radicals only.
For $C_3$ it follows by inspection that there are no integer solutions 
(or solutions with half-integer $m$) in the range considered. The same method (together with a quadrature)
applies to $C_6$, which this time yields
the solution $(n,m)=(11,10)$ in that range.

In case $C_5$ we could not find any analytical method, and the only way is the systematic search in the range $n\leq10^8$, which doesn't reveal
any integer solution. A fortiori, there seems to be no analytical methods for $k=7$ and beyond.

Of course it could be very tempting to invoke higher Number Theory to insert this problem into a general framework,
such as Shimura-Taniyama-Weil (STW) conjecture, which was used in the context of Fermat theorem, see e.g. [Da], [He].
We leave this for future investigations.

\smallskip
The paper is organized as follows: In Sect.2 we review the quadratic case of triangular numbers and compute also exact and asymptotic
solutions for the 1:st and 3:rd quartiles. In Sect.3 we investigate the case of tetrahedric numbers and beyond, up to
6-simplicial numbers, and compute the median, with the help of Mathematica. In Sect.4 we try to generalize these technics to the determination of the
1:st and 3:rd quartiles for (C3)-(C6), allowing also for $m\in{\bf N}+1/4$. 
In Sect.5 we reduce 
the equations $P_k(x,y)=0$, modulo $p=3,5,7$ when $k=3,\cdots,6$, so to save some trials when searching for integer solutions without resorting
to the algorithms set up in Sect.3.
In Sect.6 we consider the asymmetric Pascal triangle.
In Appendix we give a table of all relevant $P_k(x,y)$
obtained so far, together with the corresponding number to be approximated in the Diophantine sense.
\medskip
\noindent {\it Acknowledgements}: We thank Michel Waldschmidt and Yves Aubry for their advice. 
\medskip
\noindent {\bf 2/ Triangular numbers}
\smallskip
Eq. (Q1) actually solves simply as $m=(n+1)/2$, and Eq. (C2), (Q3) are Pell-Fermat equations. 
For (C2) (the equation considered by Ramanujan) we find $n(n+1)=2m^2$ or
$$(2n+1)^2=8m^2+1\leqno(2.1)$$
It is well-known that all ``physical'' solutions, besides this for which $n(n+1)=2m^2=0$, 
are given by the sequence indexed over $\alpha\in{\bf N}^*$
$$n_\alpha={1\over4}\bigl((3-2\sqrt2)^\alpha+(3+2\sqrt2)^\alpha-2\bigr), \quad
m_\alpha=-{\sqrt2\over8}\bigl((3-2\sqrt2)^\alpha-(3+2\sqrt2)^\alpha\bigr)\leqno(2.1)$$
Here we list the first pairs $(m,n)$
$$\eqalign{
&m_\alpha = 1; 6; 35; 204; 1189; 6930; 40\,391; 235\,416; 1372\,105; 7997\,214; \cdots\cr
&n_\alpha  = 1; 8; 49; 288; 1681; 9800; 57121; 332\,928; 1940\,449; 11\,309\,768; \cdots
}$$
They were obtained by Ramanujan by expanding $\sqrt8$ as a continued fraction (see the first approximations in Appendix). The point
is that the error term introduced in replacing $\sqrt8$ by its Diophantine approximation $(2n+1)/m$ is exactly balanced by the remainder $1/8m^2$.
This situation is actually exceptional. Here we recall the sequence of 
integers $p,q$ is a Diophantine approximation for $x\in{\bf R}^+$, iff 
$$|x-{p\over q}|<{1\over q^2}\leqno(2.2)$$
For (Q3) we find $3n(n+1)=2m(2m+1)$, and excluding the trivial solutions for which $3n(n+1)=2m(2m+1)=0$, 
$$(4m+1)^2=3(2n+1)^2-2\leqno(2.3)$$ 
This is again a Pell-Fermat Eq. Let $y=4m+1$,
$x=2n+1$.  Here we need to allow for half-integer values of $m$. Formal Calculus with Mathematica gives a fundamental sequence:
\smallskip
\noindent {\bf Proposition 2.1}: {\it Diophantine equation $y^2=3x^2-2$ has a sequence of solutions, indexed by
$\alpha\in{\bf N}^*$, of the form
$$\eqalign{
&y_\alpha={1\over2}\bigl((2+\sqrt3)^\alpha(-1+\sqrt3)-(1+\sqrt3)(2-\sqrt3)^\alpha\bigr)\cr
&x_\alpha={1\over6}\bigl((2+\sqrt3)^\alpha(3-\sqrt3)+(3+\sqrt3)(2-\sqrt3)^\alpha\bigr)
}\leqno(2.4)$$
and they are the only ones when $x<10^8$.}
\smallskip 
Here we list the 12 first couples $(m_\alpha,n_\alpha)$. Actually the $m_\alpha$'s come in consecutive pairs of integers and half-integers. 
$$\eqalign{
&m_\alpha=4.5; 17.5; 66; 247; 922.5; 3443.5; 12\,852; 47\,965; 179\,008.5; 668\,069.5; 2\,493\,270;\cr
&n_\alpha=5; 20; 76; 285; 1065; 3976; 14\,840; 55\,385; 206\,701; 771\,420; 2\,878\,980;\cr
}$$
We check that all such ${y_\alpha\over x_\alpha}$ belong to the sequence of Diophantine approximation of $\sqrt3$.
But using (2.2) we see that all approximants (not only (2.4)) give a sequence of quasi-solutions, in the sense 
$$\big|{(4m_\alpha+1)^2\over (2n_\alpha+1)^2}-3-{2\over(2n_\alpha+1)^2}\bigl|\leq{C\over n^2_\alpha}\leqno(2.5)$$
that is 
$${n_\alpha+1\choose2}-{4\over3}{m_\alpha\choose2}-m_\alpha={\cal O}(1), \ \alpha\to\infty\leqno(2.6)$$
\medskip
\noindent {\bf 3/ Tetrahedric numbers and beyond}.
\smallskip
In this Section we restrict to the median of the distribution of the $c_{nk}$, the other quartiles being investigated in Sect.4. 
\smallskip
\noindent {\it a) Tetrahedric numbers}.
\smallskip
We need to solve
$${n+2\choose3}=2{m+1\choose3}+{m+1\choose2}\leqno(C3)$$
or $n(n+1)(n+2)=m(m+1)(2m+1)$

The condition $n,m\geq2$ implies to remove the trivial non negative
solutions $(m,n)=(0,0),(1,1)$. We need also to remove negative values of $n,m$ satisfying (C3), in particular
those for which $n(n+1)(n+2)=m(m+1)(2m+1)=0$. This situation occurs for all $k$, and (Ck) will always have
at least a finite number of rational points.

Eq. (C3) can be rewritten as
$$4(n+1)^3-4(n+1)=(2m+1)^3-(2m+1) \leqno(3.2)$$
This time, Mathematica gives no sequence $(m_\alpha,n_\alpha)$, but suggests instead to solve an equation of degree 3. 
If $P(x)$ is a polynomial of degree 3, with a  
positive discriminant $4p^3+27q^2$, it is known that the equation
$P(x)=0$ has only one real root, which moreover can be expressed with real radicals of degree 2 and 3.  This is indeed the case.
\smallskip
\noindent{\bf Proposition 3.1}: {\it Consider equation $4(n+1)^3=(2m+1)^3-(2m+1)+4(n+1)$ with unknown $n\in{\bf R}$. 
For $m+1\in{\bf R}^+$, let $A=27(m+1)-81(m+1)^2+54(m+1)^3$, and $B=\sqrt{-108+729(m+1-3(m+1)^2+2(m+1)^3)^2}$. Then 
$n$ is given by }
$$n+1={2^{1/3}\over(A+B)^{1/3}}+{(A+B)^{1/3}\over 3\cdot 2^{1/3}}\leqno(3.3)$$
\smallskip
Still within the range $n\leq10^8$, it follows by inspection that there are no integer solutions. 
But Diophantine approximation of $4^{1/3}$ still gives sequences of integers $(m_\alpha,n_\alpha)$, verifying (3.2) mod ${\cal O}(n)$. 
The first terms are
$$\eqalign{
&m+1= 10; 14; 114; 391; 1903; 2407; 74\,098; \cr
&n+1= 12; 17; 143; 492; 2397; 3032; 93\,357; \cr
}$$
Look now for $m+1=p+{1\over2}$ half-integer, $p\in{\bf N}^*$, we rewrite (3.2) in the form
$$4n(n+1)(n+2)=(2p-1)2p(2p+1)\leqno(3.4)$$
Let $q=(2p-1)2p(2p+1)$. Inserting into (3.3) we find
$$n+1={4+3\bigl(q+\sqrt{-{64\over27}+q^2}\bigr)^{2/3}\over 6\bigl(q+\sqrt{-{64\over27}+q^2}\bigr)^{1/3}}\leqno(3.5)$$
One checks that (3.5) holds for any $n\geq1$, with $q=q(n)=4n(n+1)(n+2)$. 
So we just recover (3.4), which is neither fulfilled for integer $p$, in the range $n\leq10^8$.
So there are no solution to (C3) with $m$ half-integer.

Still again, Diophantine approximations of $2^{1/3}$ give sequences of integers $(m_\alpha,n_\alpha)$, verifying (3.4) mod ${\cal O}(n)$. 
The first terms are
$$\eqalign{
&m+1= 3.5; 4.5; 23.5; 27.5; 50.5; 227.5;\cr
&n+1= 4; 5; 29; 34; 63; 286;
}$$
\smallskip
\noindent {\it b) 4-simplicial numbers}.
\smallskip
Consider 4-simplicial numbers $C_4$, i.e.
$${n+3\choose4}=2{m+2\choose4}+{m+2\choose3}\leqno(C4)$$
which we rewrite as
$(n+3)(n+2)(n+1)n=2(m+2)(m+1)^2m$,
The condition $n,m\geq2$ implies again to remove the trivial non negative
solutions $(m,n)=(0,0), (1,1)$.

This simplifies to
$$(2n+3)^4-32(m+1)^4+32(m+1)^2-10(2n+3)^2+9=0\leqno(3.6)$$
a quadratic equation in 
$i=(2n+3)^2\geq25, j=(m+1)^2\geq9$, so we start to solve Pell-Fermat equation
$(i-5)^2=32j^2-32j+16$, or with $u=i-5$
$$u^2=8(2j-1)^2+8\leqno(3.7)$$
from which we remove the solutions $(u,j)=(4,1), \ (20,4)$. 
Mathematica gives us four fundamental sequences.
\smallskip
\noindent {\bf Proposition 3.2}: Diophantine Eq. $u^2=8(2j-1)^2+8$ has at least 4 sequences of integer solutions indexed by 
$\alpha\in{\bf N}^*$, namely
$$\leqalignno{
&u_\alpha=(577-408\sqrt2)^\alpha(2+\sqrt2)-(-2+\sqrt2)(577+408\sqrt2)^\alpha\cr
&j_\alpha={1\over4}\bigl(2+(-1+\sqrt2)(577+408\sqrt2)^\alpha-(577-408\sqrt2)^\alpha(1+\sqrt2)\bigr)\cr
&(u,j)=(4,1), \ (20,4), \ (676,120), \ (780\,100,137\,904), \ \&c &(3.8)\cr
&u_\alpha=(577-408\sqrt2)^\alpha(10+7\sqrt2)+(10-7\sqrt2)(577+408\sqrt2)^\alpha\cr
&j_\alpha={1\over4}\bigl(2-(7+5\sqrt2)(577-408\sqrt2)^\alpha-(577+408\sqrt2)^\alpha(7-5\sqrt2)\bigr)\cr
&(u,j)=(4,1), \ (20,4), \ (116,21), \ (133\,844, 23\,661), \ \&c &(3.9)\cr
&u_\alpha=(577+408\sqrt2)^\alpha(2+\sqrt2)-(-2+\sqrt2)(577-408\sqrt2)^\alpha\cr
&j_\alpha={1\over4}\bigl(2+(1+\sqrt2)(577+408\sqrt2)^\alpha-(577-408\sqrt2)^\alpha(-1+\sqrt2)\bigr)\cr
&(u,j)=(4,1), \ (20,4), \ (3\,940,697), \ (4\,546\,756, 803\,761), \ \&c &(3.10)\cr
}$$
$$\leqalignno{
&u_\alpha=(577+408\sqrt2)^\alpha(10+7\sqrt2)-(-10+7\sqrt2)(577-408\sqrt2)^\alpha\cr
&j_\alpha={1\over4}\bigl(2+(7+5\sqrt2)(577+408\sqrt2)^\alpha-(577-408\sqrt2)^\alpha(-7+5\sqrt2)\bigr)\cr
&(u,j)=(4,1), \ (20,4), \ (22\,964,460), \ (26\,500\,436, 4\,684\,660), \ \&c &(3.11)
}$$
Except for $(u,j)=(4,1), \ (20,4)$ we have excluded, it seems by inspection that $(u+5,j)$ are never perfect squares. 
We can consider also consider the equation generalizing (3.7)
$$u^2=8v^2+8\leqno(3.12)$$ 
which admits also a fundamental sequence of integer solutions
$$\eqalign{
&u_\alpha=(3-2\sqrt2)^\alpha(2+\sqrt2)-(-2+\sqrt2)(3+2\sqrt2)^\alpha\cr
&v_\alpha={1\over2}\bigl((3+2\sqrt2)^\alpha(-1+\sqrt2)-(1+\sqrt2)(3-2\sqrt2)^\alpha\bigr)
}\leqno(3.13)$$
but neither leads to any integer solution to (3.7) (up to $n\leq10^8$). 
So we try $v=w-1/2$ half-integer, which gives $2j=w+1/2$. Condition $j=(m-1)^2$ for $m$ half-integer, $m=p+1/2$, implies
$w+1/2=2(p-1/2)^2=2p^2-2p+1/2$, this gives the quadratic Eq. $2p^2-2p-w=0$, which has an integer solution iff $2w+1$ is a perfect square.
But when $v=w-1/2$, (3.12) can be written as
$u^2=8w^2-8w+10$, or else
$$u^2=2(2w-1)^2+8\leqno(3.14)$$
which is again a Pell-Fermat Eq. But contrary to (3.7) or (3.12) there are no integer solution to (3.14), cf. [Chr,p.483]. 

Still again, Diophantine approximation gives sequences of integers or half-integers $(m_\alpha,n_\alpha)$, verifying (C4) mod ${\cal O}(n^2)$.
\medskip
\noindent {\it c) 5-simplicial numbers}.
\smallskip
Consider 5-simplicial numbers $C_5$, i.e.
$${n+4\choose5}=2{m+3\choose5}+{m+3\choose4}\leqno(C5)$$
with the condition $n,m\geq2$, which we rewrite as
$(n+4)(n+3)(n+2)(n+1)n=(m+3)(m+2)(m+1)m(2m+3)$,
and we remove again the non negative solutions $(m,n)=(0,0),(1,1)$. This we rewrite as
$P_5(n+2,2m+3)=0$ where $P_5$ as in (A.2). 
Mathematica gives no hint at solving this 5:th degree equation, and systematic search up to $n\leq10^8$ either gives no integer solutions.  
But Diophantine approximation gives sequences of integers or half-integers $(m_\alpha,n_\alpha)$, verifying (C5) mod ${\cal O}(n^3)$.
\medskip
\noindent {\it d) 6-simplicial numbers}
\smallskip
Consider now 6-simplicial numbers, i.e.
$${n+5\choose6}=2{m+4\choose6}+{m+4\choose5}\leqno(C6)$$
with the condition $n,m\geq2$. This can be rewritten as
$(n+5)(n+4)(n+3)(n+2)(n+1)n=2(m+4)(m+3)(m+2)^2(m+1)m$, and we remove the solutions $(m,n)=(0,0),(1,1)$.
The methods elaborated for (C3) will give at least 2 non-trivial solutions
in the range $n\leq10^8$. Let 
$x=2n+5$, $y=2m+4$, (C6) leads to
$$x^6-2y^6-35x^4+40y^4+259x^2-128y^2-225=0$$
or if we let $u=x^2$, $v=y^2$
$$u^3-2v^3-35u^2+40v^2+259u-128v-225=0\leqno(3.17)$$
As in Proposition 3.1 for (C3)  Mathematica solves (3.17) by real radicals.
\smallskip
\noindent{\bf Proposition 3.3} {\it For any $u\in{\bf R}_+$, let
$$A=5770+13\,986u-1890u^2+54u^3, \  B=2\sqrt{-143\,982\,592+\bigl(2885+27u(259+(-35+u)u\bigr)^2}$$ 
Then 
$u^3-2v^3-35u^2+40v^2+259u-128v-225=0$ has the unique solution
$v\in{\bf R}_+$ given by }
$$v={20\over3}+{416\over3(A+B)^{1/3}}+{1\over6}(A+B)^{1/3}\leqno(3.18)$$
\smallskip
For $v\leq 10^8$, it follows by inspection of (3.18), that the only integer solutions of (3.17), 
which are also perfect squares, are $(u,v)=(49, 36), (729, 576)$.
They correspond to $(n,m)=(1,1)$ (which is excluded) and $(n,m)=(11,10)$. 
There are no solution of (3.17) in this interval with $m$ half-integer. 
Again Diophantine approximation gives sequences of integers or half-integers $(m_\alpha,n_\alpha)$, verifying (C6) mod ${\cal O}(n^4)$.
\medskip
\noindent {\bf 4/ First and third quartiles for distributions of the 3,4,5,6-simplicial numbers}.
\smallskip
Equations (Q1) and (Q3) generalize to higher orders, leading to generalized Pell-Fermat equations with additional lower order terms 
provided we allow also for $m\in{\bf N}+1/4$. Even if there are no exact integer solution, we can still find approximations in the Diophantine sense.

For $C_3$ the 1:st quartile is defined by
$${n+2\choose3}=4{m+1\choose3}+{m+1\choose2}\leqno(Q13)$$
with the condition $n,m\geq2$.  This can be rewritten as
$(n+1)(n+2)n=(4m-1)(m+1)m$, and we remove again the solutions $(m,n)=(0,0), (1,1)$.
With $m+1=k+1/4$, the latter equation takes the form
$$8(n+1)^3=4(2k-1)^3+8(n+1)-7(2k-1)+3\leqno(4.1)$$
which is a generalized Pell-Fermat equation with a constant term. 

The 3:rd quartile is defined by
$${n+2\choose3}={4\over3}{m+1\choose3}+{m+1\choose2}\leqno(Q33)$$
with the condition $n,m\geq2$. With $m+1=k-1/4$, a similar computation leads to
$$24(n+1)^3=4(2k-1)^3+24(n+1)-7(2k-1)-3\leqno(4.2)$$
which is a generalized Pell-Fermat equation with a constant term.
\smallskip
\noindent{\bf Proposition 4.1}: {\it Consider the quartiles (Q13) and (Q33) characterized by (4.1) and (4.2) resp. Then as in Proposition 3.1,
the 3:rd degree equations can be solved in the following form: (\dots)}

\vskip 2truecm

For $C_4$ the 1:st quartile is defined by
$${n+3\choose4}=4{m+2\choose4}+{m+2\choose3}\leqno(Q14)$$
with the condition $n,m\geq2$. Reducing to a generalized Pell-Fermat equation with lower order terms
leads to
$$(2n+3)^4-10(2n+3)^2+9=4(2m+{3\over2})^4+106(2m+{3\over2})^2+2772(2m+{3\over2})+23409/4$$
Trying $m+2=k-{1\over4}$, $k\in{\bf N}$, we get
$$(2n+3)^4-10(2n+3)^2=4(2k-3)^4+106(2k-3)^2+2772(2k-3)+23409/4\leqno(4.10)$$
which clearly has no integer solution, so is the case
$m+2=k+{1\over4}$.

The 3:rd quartile is defined by
$${n+3\choose4}={4\over3}{m+2\choose4}+{m+2\choose3}\leqno(Q34)$$
with the condition $n,m\geq2$. Reducing to a generalized Pell-Fermat equation with lower order terms leads to
$$3\bigl[(2n+3)^4-10(2n+3)^2+9\bigr]=4\bigl[(2m+{5\over2})^4-{11\over2}(2m+{5\over2})^2-3(2m+{5\over2})-{27\over8}\bigr]$$
As above the substitution $m+2=k+{1\over4}$, $k\in{\bf N}$ gives
$$6(2n+3)^4-60(2n+3)^2=8(2k-1)^4-44(2k-1)^2-24(2k-1)-81\leqno(4.11)$$
which clearly has no integer solution, so is the case
$m+2=k-{1\over4}$.

For $C_6$ the 3:rd quartile is defined by
$${n+5\choose6}={4\over3}{m+4\choose6}+{m+4\choose5}\leqno(Q36)$$
with the condition $n,m\geq2$. We proceed as before but the computation is somewhat more tedious. We find
$$\eqalign{
3&\bigl[(2n+5)^6-35(2n+5)^4+259(2n+5)^2-225\bigr]=\cr
&4\bigl[(2m+8)^6-21(2m+8)^5+160(2m+8)^4-540(2m+8)^3+784(2m+8)^2-384(2m+8)\bigr]
}$$
To eliminate the 2:nd term on the RHS we expand $(2m+8-\alpha)^6$ and are led to choose $\alpha=7/2$. With $m+8=k\pm{1\over4}$
we find $2m+8-7/2=2k-3$, resp. $2m+8-7/2=2(k-2)$.
When $m+8=k+{1\over4}$, we substitute $2k-3$ to $2m+8$ several times and find eventually
$$\eqalign{
3&(2n+5)^6-105(2n+5)^4+777(2n+5)^2=\cr
&4(2k-3)^6-95(2k-3)^4-60(2k-3)^3+{1939\over4}(2k-3)^2+{1500\over4}(2k-3)+{6075\over16}
}\leqno(4.12)$$
which is a generalized Pell-Fermat equation with lower order terms, but without a integer solution.
The same conclusion holds for $m+8=k-{1\over4}$.
However there are always asymptotic solutions in the Diophantine sense, for instance for (4.12) it suffices to compute the
approximants of $(4/3)^{1/6}$ as a continued fraction.

Note that it could turn out that there could still be integer solutions to (Q36) but not in the form $m+8=k\pm{1\over4}$
(which provides asymptotic solutions).

Consider at last the 1:st quartile
$${n+5\choose6}=4{m+4\choose6}+{m+4\choose5}\leqno(Q16)$$
\medskip
\noindent {\bf 5/ Some arithmetics}.
\smallskip
A straighforward way to find solutions to (Ck) is by trial. Using Mathematica we can easily check all values of $n$ less than
$10^8$ say. None of them is a solution in that range, except the one for (C6).
However to restrain somewhat the domain of numbers $n,m$ we are looking for, we can first
reduce Eq. (Ck) modulo prime numbers $p$, using the simple fact that if $a\in{\bf N}^*$,
then $a^p=a[p]$, and if moreover $p$ does not divide $a$, then $a^{p-1}=1[p]$. This reduces sometimes up to 50\% the range of trials.

Apply this to (C3) (tetrahedral numbers) with $p=3$, we let $x=n+1, y=2m+1$ so that (3.2) reads $4(x^3-x)=y^3-y$, which always holds
in ${\bf F}_3$, so this gives no criterion at all.

Consider now (C4) with $p=5$. With $x=2n+3,y=m+1$, we recall from (3.6) the identity
$x^4-10x^2=32y^4-32y^2-9$. When 5 does not divide $x$ and $y$, we find $32y^2-10x^2=2$ in ${\bf F}_5$. By inspection we find 
that this equation has no solution in ${\bf F}_3$ when
$(x,y)\in \{1,2,3,4\}\times\{2,3\}$. When $x$ or $y$ vanish in ${\bf F}_3$, we find either that there are no solutions such that
$(x,y)\in\{(0,1),(0,4),(0,0)\}$. So we can save about 50\% of trials.

Consider (C5) with $p=5$. With $x=n+2,y=2m+3$, we recall from (3.15) the identity
$16x^5-80x^3+64x=y^5-10y^3+9y$. Reducing mod 5, we find $80(x^3-x)=10(y^3-y)$, which always holds
in ${\bf F}_5$, so again this gives no criterion.

Consider (C6) with $p=7$. With $x=2n+5, y=2m+4$, recall from (3.16) the identity
$x^6-35x^4+259x^2=2y^6-40y^4+128y^2+225$. When 7 does not divide $x$ and $y$, we find $40y^4-128y^2-226=0$ in ${\bf F}_7$.
By inspection we find that this equation has no solution when $(x,y)\in{\bf F}_7^*\times\{1,6\}$.
When 7 divides $y$ but not $x$, the equation reduces to $x^6-1=0$, which always holds.
When 7 divides $x$, the equation reduces to
$2y^6-40y^4+128y^2+225=0$ in ${\bf F}_7$
we find that there are no solutions whenever $y\in\{2,3,4,5\}$.
So again we can save about 50\% of trials. Note that (C6) reduces to a polynomial equation of only one variable in ${\bf F}_7$,
which may give a hint at the special case played by the 6-simplicial numbers.

\medskip
\noindent {\bf 6/ Asymmetric Pascal triangle}
\smallskip
We consider floors $C_2^a$ and $C_3^a$ only, and content with Diophantine approximation, using also that, as $n$ becomes large, 
the leading term of ${n\choose k}_a$ is ${n\choose k}$. 

For $C_2^a$ consider
${n+1\choose 2}_a={n+1\choose 2}+an$, which we take to be equal to
$2{m\choose 2}_b+m+b=2{m\choose 2}+2b(m-1)+m+b$ mod ${\cal O}(1)$ as $n\to\infty$. This gives
$${n+1\choose 2}-2{m\choose 2}-m=2bm-an-b+{\cal O}(1)$$
We know that when ${2m+1\over m}$ approximates $\sqrt8$ at this order, the LHS is 0, so we may replace $m$ by $2n\over\sqrt8$ in the RHS
which gives $\bigl({4\over\sqrt8}b-a\bigr)n+{\cal O}(1)$, so we choose $b=a2^{-1/2}$.

For $C_3^a$ consider
${n+1\choose 3}_a={n+1\choose 3}+a{n\choose 2}$, which we take to be equal mod ${\cal O}(n)$ to
$$2{m\choose 3}_b+{m\choose 2}_b=2{m\choose 3}+2b{m-1\choose 2}+{m\choose 2}+b(m-1)$$
Using again Diophantine approximation of $4^{1/3}$ by ${2m-1\over n}$ to this order, we need as before to cancel the term
$bm^2-an^2/2$, which gives $b=a2^{-1/3}$. So we proved
\smallskip
\noindent {\bf Proposition 5.1}: At floor $C_2^a$, for $b=a2^{-1/2}$, we have
$${n+1\choose 2}_a=2{m\choose 2}_b+m+b+{\cal O}(1)$$
for the sequence $(m_\alpha,n_\alpha)$ given in (2.2), 
while for $b=a2^{-1/3}$ at floor $C_3^a$, we have
$${n+1\choose 3}_a=2{m\choose 3}_b+{m\choose 2}_b+{\cal O}(n)$$
for the sequence given by Diophantine approximation (this holds for integer and half-integer $m$). 
\medskip
\noindent {\bf Appendix}
\smallskip
\noindent {\it a) Table of Diophantine approximations}.
$$\eqalign{
&2^{3/2}=3-{1\over 6-}\,{1\over6-}\,{1\over6-}\,\cdots=3; \ {17\over6}; \ {99\over35}; \ {577\over204}; \cr
&2^{2/3}=1+{1\over 1+}\,{1\over1+}\,{1\over2+}\,{1\over2+}\,\cdots= 1; {3\over2}; {8\over5}; {19\over22}; {27\over17}; {100\over63}; 
{227\over143}; {781\over492}; {1008\over635}; {3805\over2397}; {4813\over3032}; {14\,8195\over93\,357};\cr
&2^{1/3}=1; {4\over3}; {5\over4}; {29\over23}; {34\over27}; {63\over50}; {286\over227}; \cr
&2^{1/2}=1+{1\over 2+}\,{1\over2+}\,{1\over2+}\,\cdots; \ {17\over12}; \ {99\over70}; \ {577\over408}; \cr
&(4/3)^{1/6}=
}\leqno(A.1)$$
\smallskip
\noindent {\it b) Some plane algebraic curves}. 
\smallskip
We have met the following polynomials, for which we can always find approximate integer roots in the Diophantine sense, when we compute
approximants of the irrational number to the right:
$$\matrix{
P_2(x,y)=x^2-8y^2-1&x=2n+1&y=m&2^{3/2}\cr
Q_2(x,y)=x^2-3y^2+2&x=2n+1&y=4m+1&3^{1/2}\cr
P_3(x,y)=4x^3-y^3-4x+y&x=n&y=2m-1&2^{2/3}\cr
P_4(x,y)=x^4-32y^4+32y^2-10x^2+9&x=2n-1&y=m-1&2^{5/4}\cr
P_5(x,y)=16x^5-y^5+10y^3-80x^3+64x-9y&x=n-1&y=2m-3&2^{4/5}\cr
P_6(x,y)=x^6-2y^6-35x^4+40y^4+259x^2-128y^2-225&x=2n-3&y=2m-4&2^{1/6}\cr
}\leqno(A.2)$$
\medskip
\noindent{\bf References}:

\noindent [Be] A.Beiler. ``The Pellian.'' Ch. 22 in Recreations in the Theory of Numbers: 
The Queen of Mathematics Entertains. New York: Dover, pp. 248-268, 1966.

\noindent [Ch] G.Chrystal. Textbook of Algebra, 2nd ed., Vol. 2. New York: Chelsea, pp. 478-486, 1961.

\noindent [Da] H.Darmon. A Proof of the Full Shimura-Taniyama-Weil Conjecture Is Announced. Notices AMS, p.1397-1401, Dec.1999.

\noindent [Di] L.Dickson. ``Pell Equation: Made Square.'' Ch. 12 in History of the Theory of Numbers, 
Vol. 2: Diophantine Analysis. New York: Dover, pp. 341-400, 2005.

\noindent [GelKaZe] I.M.Gelfand, M.Kapranov, A.Zelevinsky.
Discriminants, Resultants, and Multidimensional Determinants. Mathematics: Theory \& Applications. Birkh\"auser, Boston, MA, 1994

\noindent [Ge] A.G\'erardin. Formules de r\'ecurrence. Sphinx-Oedipe 5, p.17-29, 1910. 

\noindent [He] Y.Hellegouarch. Invitation aux math\'ematiques de Fermat-Wiles. Masson, Paris, 1997.

\noindent [Ran] B.Rand\'e. Les carnets indiens de Srinivasa Ramanujan. Cassini, Paris 2002.

\noindent [Wa] M.Waldschmidt. Pell's Equation.
https://webusers.imj-prg.fr/~michel.waldschmidt/articles/ pdf/BamakoPell2010.pdf

\noindent [Wo] Wolfram. Pell Equation. mathworld.wolfram.com/PellEquation.html

\bye

\bye

%% file: Definitions.tex
\magnification=1100

\hsize 17truecm
\vsize 22truecm

\font\twelvec=msbm10 at 12pt
\font\sevenc=msbm10 at 9pt
\font\fivec=msbm10 at 7pt

\newfam\co
\textfont\co=\twelvec
\scriptfont\co=\sevenc
\scriptscriptfont\co=\fivec

\def\Sum{\displaystyle\sum}

\baselineskip 15pt
